\documentclass{amsart}
\usepackage{amsfonts}
\usepackage{amsmath}
\usepackage{enumerate}
\usepackage{amsthm}
\usepackage{hyperref}
\usepackage{cite}

\begin{document}

\newtheorem{definition}{Definition}[section]
\newtheorem{theorem}[definition]{Theorem}
\newtheorem{proposition}[definition]{Proposition}
\newtheorem{remark}[definition]{Remark}
\newtheorem{lemma}[definition]{Lemma}
\newtheorem{corollary}[definition]{Corollary}
\newtheorem{example}[definition]{Example}

\numberwithin{equation}{section}

\title[Prescribed Chern scalar curvatures]{Prescribed Chern scalar curvatures on compact Hermitian manifolds with negative Gauduchon degree}
\author{Weike Yu}

\date{}

\begin{abstract}
In this paper, we investigate the problem of prescribing Chern scalar curvatures on compact Hermitian manifolds with negative Gauduchon degree. By studying the convergence of the associated geometric flow, we obtain some existence results when the candidate curvature function is nonzero and nonpositive. Furthermore, we also consider the case that the candidate curvature function is sign-changing, and establish some existence and nonexistence results.
\end{abstract}
\keywords{Prescribed Chern scalar curvature problem;  Hermitian manifold; Chern-Yamabe flow; Trudinger-Moser inequity.}
\subjclass[2010]{53B35, 53C44, 32Q99.}
\maketitle
\section{Introduction}
Suppose that $(M^n, J, h)$ is a compact complex manifold of complex dimension $n\geq2$ with a Hermitian metric $h$. The fundamental form $\omega$ of the Hermitian metric $h$ is given by $\omega(\cdot, \cdot)=h(J\cdot, \cdot)$. In this paper, we will confuse the Hermitian metric $h$ and its corresponding fundamental form $\omega$. For a Hermitian metric $\omega$, there exists a $1$-form $\theta$ on $M$ (called the Lee form or torsion $1$-form) such that
\begin{align}
d\omega^{n-1}=\theta\wedge \omega^{n-1}.
\end{align}
A Hermitian metric $\omega$ is said to be balanced if $\theta=0$. It is said to be Gauduchon if $d^*\theta=0$, where $d^*$ denotes the adjoint operator of $d$ with respect to $\omega$.

On a Hermitian manifold $(M^n, \omega)$, there is a unique linear connection $\nabla^{Ch}$ preserving the Hermitian metric and the complex structure, whose torsion $T^{Ch}$ has vanishing $(1,1)$-part everywhere. Such a connection is called the Chern connection. Then, the scalar curvature with respect to the Chern connection (referred to as the Chern scalar curvature) can be written as
\begin{align}
S^{Ch}(\omega)=\text{tr}_{\omega}\sqrt{-1}\bar{\partial}\partial \log{\omega^n}, 
\end{align}
where $\omega^n$ denotes the corresponding volume form.

On a compact Hermitian manifold $(M^n, \omega)$, the following problem is a Hermitian analogue of prescribing scalar curvature: Given a smooth real-valued function $g$ on $M$, does there exist a Hermitian metric $\hat{\omega}$ conformal to $\omega$, that is, $\hat{\omega}=e^{\frac{2}{n}u}\omega$, such that its Chern scalar curvature $S^{Ch}(\hat{\omega})=g$? According to \cite{[ACS]}, we know that it is equivalent to solving the following partial differential equation:
\begin{align}
-\Delta^{Ch}_\omega u+S^{Ch}(\omega)=ge^{\frac{2}{n}u},\label{1.3}
\end{align}
where $\Delta^{Ch}_\omega$ denotes the Chern Laplacian with respect to $\omega$. When $g$ is a constant, the above problem is referred to as the Chern-Yamabe problem, which was first proposed by Angella-Calamai-Spotti in \cite{[ACS]}, and they proved that it is solvable on any compact Hermitian manifolds with the Gauduchon degree $\Gamma(\{\omega\})\leq 0$ (see \eqref{2.6.}). In \cite{[LY]}, by using the method of upper and lower solutions, Liu-Yao showed that an equation of the type like \eqref{1.3} has a smooth solution when $S(\omega)$ is a negative constant and $g$ is nonpositive and not equal to zero identically on $M$ (denoted by $g\leq 0\ (\not\equiv0)$). Hence, combining \cite[Theorem 4.1]{[ACS]}, we know that every smooth function $g\leq 0\ (\not\equiv 0)$ can be prescribed as the Chern scalar curvature in the conformal class, provided that $\Gamma(\{\omega\})< 0$ (also see \cite[Theorem 2.5]{[Fus]}). Moreover, Fusi \cite{[Fus]} also obtained that if $\omega$ is a balanced metric with zero Chern scalar curvature, then every smooth sign-changing function $g$ with $\int_M g\frac{\omega^n}{n!}<0$ is the Chern scalar curvature of a metric conformal to $\omega$.

In \cite{[ACS]}, the authors proposed to use the following flow (called the Chern-Yamabe flow) to tackle the above Chern-Yamabe problem:
\begin{equation}
\left\{
\begin{aligned}
&\frac{\partial u}{\partial t}=\Delta^{Ch}_\omega u-S(\omega)+\lambda e^{\frac{2}{n}u},\\
&u(x,0)=0,
\end{aligned}
\right.
\end{equation}
where $\lambda$ is a constant. Later, by using the comparison principle, Lejmi-Maalaoui \cite{[LM]} proved that when the Gauduchon degree is negative, the above flow converges to a solution of the Chern-Yamabe problem, which recovered the result of Angella-Calamai-Spotti. In \cite{[CZ]}, Calamai-Zou introduced a different flow to study the Chern-Yamabe problem. By using geometric flows related to Calamai-Zou's Chern–Yamabe flow, Ho \cite{[Ho]} studied the problem of prescribing Chern scalar curvature on balanced Hermitian manifolds with negative Chern scalar curvature. Besides, Ho-Shin \cite{[HS]} showed that the solution to the Chern-Yamabe problem is unique under suitable conditions and obtained some results related to the Chern-Yamabe soliton.

In this paper, we will focus on the prescribed Chern scalar curvature problem on compact Hermitian manifolds with negative Gauduchon degree. Firstly, we will investigate the convergence of the following geometric flow:
\begin{equation}
\left\{
\begin{aligned}
&\frac{\partial u}{\partial t}=\Delta^{Ch}_\eta u-s_0+ge^{\frac{2}{n}u},\\
&u(x,0)=0,
\end{aligned}
\right.
\label{2.1}
\end{equation}
where $\eta$ is the unique Gauduchon representative in the conformal class $\{\omega\}$ with volume $1$ (cf. Theorem \ref{theorem2.2}), $s_0$ is the Chern scalar curvature with respect to $\eta$, and $g$ is a smooth function on $M$ with $g\leq 0\ (\not\equiv0)$. More precisely, we prove that
 \begin{theorem}\label{theorem1.1}
Let $(M^n,\omega)$ be a compact Hermitian manifold with complex dimension $n\geq 2$ and Gauduchon degree $\Gamma (\{\omega\})<0$, and assume that $g\in C^\infty(M)$ is nonpositive and not equal to zero identically. Then the flow \eqref{2.1} has a unique solution $u$ defined on $M\times [0,+\infty)$. Furthermore, Suppose either
\begin{enumerate}
\item $g$ is negative on $M$
\end{enumerate}
or
\begin{enumerate}[\ \ \ \ (2)]
\item $g$ is nonpositive and not equal to zero identically on $M$ and $\{\omega\}$ contains a balanced metric,
\end{enumerate}
then there is a subsequence $t_i\rightarrow +\infty$ such that $u(\cdot,t_i)$ converges to $u_\infty$ in $C^{1,\alpha}(M)$, where $\alpha\in (0,1)$, and $u_\infty$ is the unique smooth solution of 
\begin{align}\label{1.6.}
-\Delta^{Ch}_\eta u_\infty+s_0=ge^{\frac{2}{n}u_\infty}.
\end{align}
\end{theorem}

Let $PC(\omega)$ denote the set of $C^\infty(M)$ functions which are Chern scalar curvatures of all $\hat\omega \in \{\omega\}$. In other words, $PC(\omega)$ is the set of smooth functions $g$ for which one can find a solution of \eqref{1.3}. Then, from Theorem \ref{theorem1.1}, we have the following corollary which partially recovers \cite[Proposition 2.1]{[LY]} or \cite[Theorem 2.5]{[Fus]}. 
\begin{corollary}\label{corollary1.3}
Let $(M^n,\omega)$ be a compact Hermitian manifold with complex dimension $n\geq 2$ and Gauduchon degree $\Gamma (\{\omega\})<0$. Suppose that $g$ is a smooth function on $M$. If either
\begin{align*}
(1)\ g<0\ \ \ \ \  \text{or}\ \ \ \ \ (2)\ g\leq 0\ (\not\equiv0)\ \text{and}\ \{\omega\}\ \text{contains a balanced metric},
\end{align*}
then $g\in PC(\omega)$.
\end{corollary}

\begin{remark}
In \cite{[Ho]}, using a different geometric flow, Ho proved that when $(M^n,\omega_0)$ is a balanced manifold with Chern scalar curvature $S^{Ch}(\omega_0)<0$, then $\{g\in C^\infty(M): g<0\}\subset PC(\omega_0)$.
\end{remark}

Now we turn to the case that the candidate curvature function $g$ is sign-changing. Since every sign-changing function $g\in C^{\infty}(M)$ can be expressed as $g=g_0+\lambda$, where $g_0$ is a nonconstant smooth function satisfying $\max_M g_0=0$ and $\lambda>0$ is a constant, then by the implicit function theorem, we have the following
\begin{theorem}\label{theorem1.4}
Let $(M^n,\omega)$ be a compact Hermitian manifold with complex dimension $n\geq 2$ and Gauduchon degree $\Gamma (\{\omega\})<0$. Let $g_0$ be a nonconstant smooth function on $M$ with $\max_M g_0=0$. Then there exists a constant $\lambda^*\in (0,-\min_M g_0]$ such that 
\begin{enumerate}
\item If $0<\lambda<\lambda^*$, then $g_0+\lambda\in PC(\omega)$.
\item If $\lambda>\lambda^*$, then $g_0+\lambda\not\in PC(\omega)$.
\end{enumerate}
\end{theorem}

We remark that whether $g_0+\lambda^*$ belongs to $PC(\omega)$ is still open for arbitrary Hermitian manifold $(M^n,\omega)$. However, when $\dim_{\mathbb{C}}M=2$, we prove that
\begin{theorem}
Let $(M^n,\omega)$ be a compact Hermitian manifold with complex dimension $n=2$ and Gauduchon degree $\Gamma (\{\omega\})<0$, and suppose that $\{\omega\}$ contains a balanced metric.  Let $g_0$ and $\lambda^*$ be as in Theorem \ref{theorem1.4}. Then $g_0+\lambda^*\in PC(\omega)$. 
\end{theorem}
\begin{remark}
In \cite{[DL]}, Ding-Liu established a similar existence result on compact Riemannian surface (i.e., $\dim_{\mathbb{C}}=1$). Unfortunately, Ding-Liu's argument is not applicable in higher dimensional cases directly. In the present paper, combining the spirt of \cite{[DL]} and the method of Moser iteration, we can prove the above theorem.
\end{remark}

\textbf{Acknowledgements.} The author would like to thank Prof. Yuxin Dong and Prof. Xi Zhang for their useful discussions and helpful comments.

\section{Preliminaries}
Let $(M^n, J, h)$ be a Hermitian manifold with $\dim_{\mathbb{C}} M=n$ and its corresponding fundamental form $\omega(\cdot,\cdot)=h(J\cdot, \cdot)$. Let $TM^{\mathbb{C}}=TM\otimes \mathbb{C}$ be the complexified tangent space $TM$, and we extend $J$ and $h$ from $TM$ to $TM^{\mathbb{C}}$ by $\mathbb{C}$-linearity. Then we have the following decomposition
\begin{align}
TM^{\mathbb{C}}=T^{1,0}M\oplus T^{0,1}M,
\end{align}
where $T^{1,0}M$ and $T^{0,1}M$ are the eigenspaces of complex structure $J$ corresponding to the eigenvalues $\sqrt{-1}$ and $-\sqrt{-1}$, respectively. Note that every $m$-form can also be decomposed into $(p, q)$-forms for each $p, q \geq 0$ with $p + q = m$ by extending the complex structure $J$ to forms.

On a Hermitian manifold $(M^n, J, h)$, there exists a unique affine connection $\nabla^{Ch}$ preserving both the Hermitian metric $h$ and the complex structure $J$, that is, $\nabla^{Ch}h=0$, $\nabla^{Ch}J=0$, whose torsion $T^{Ch}(X,Y)=\nabla_XY-\nabla_YX-[X,Y]$ satisfies $T^{Ch}(JX,Y)=T^{Ch}(X,JY)$. Note that we will confuse the Hermitian metric $h$ and its corresponding fundamental form $\omega$ in this paper. 

Given a Hermitian manifold $(M^n, \omega)$ with Chern connection $\nabla^{Ch}$, it is well-known that the corresponding Chern scalar curvature is given by
\begin{align}
S^{Ch}(\omega)=\text{tr}_{\omega}Ric^{(1)}(\omega)=\text{tr}_{\omega}\sqrt{-1}\bar{\partial}\partial \log{\omega^n},
\end{align}
where $Ric^{(1)}(\omega)$ is the first Chern Ricci curvature.

Analogous to the Laplace operator in Riemannian geometry, there is a canonical elliptic operator in Hermitian geometry, which is called the Chern Laplace operator. More precisely, for any smooth function $u: M\rightarrow \mathbb{R}$, the Chern Laplacian $\Delta^{Ch}$ of $u$ is defined by
\begin{align}
\Delta^{Ch}u=-2\sqrt{-1} tr_{\omega}\overline{\partial}\partial u.
\end{align}
\begin{lemma}[cf. \cite{[Gau]}]\label{lemma2.1}
On a compact Hermitian manifold $(M^n,\omega)$, we have
\begin{align}
-\Delta^{Ch}u=-\Delta_d u+(du,\theta)_{\omega},
\end{align}
where $\Delta_d$ is the Hodge-de Rham Laplacian, $\theta$ is the Lee form or torsion $1$-form given by $d\omega^{n-1}=\theta\wedge \omega^{n-1}$, and $(\cdot, \cdot)_\omega$ denotes the inner product on $1$-form induced by $\omega$. In particular, when $\omega$ is balanced (i.e., $\theta\equiv0$), the Chern Laplacian and the Hodge-de Rham Laplacian on smooth functions coincide.
\end{lemma}
Let 
\begin{align}
\{\omega\}=\{e^{\frac{2}{n}u}\omega\ |\ u\in C^\infty(M)\}
\end{align}
denote the conformal class of the Hermitian metric $\omega$. In \cite{[Gau1]}, Gauduchon proved the following interesting result
\begin{theorem}\label{theorem2.2}
Let $(M, \omega)$ be a compact Hermitian manifold with $\dim_{\mathbb{C}} M\geq 2$, then there exists a unique Gauduchon metric $\eta\in \{\omega\}$ (i.e., $d^*\theta=0$) with volume $1$.
\end{theorem}
In terms of the above theorem, one can define an invariant $\Gamma(\{\omega\})$ of the conformal class $\{\omega\}$ as follows:
\begin{align}\label{2.6.}
\Gamma(\{\omega\})=\frac{1}{(n-1)!}\int_M c^{BC}_1(K^{-1}_M)\wedge\eta^{n-1}=\int_M S^{Ch}(\eta)d\mu_\eta,
\end{align}
where $c^{BC}_1(K^{-1}_M)$ is the first Bott-Chern class of anti-canonical line bundle $K^{-1}_M$, and $d\mu_\eta$ denotes the volume form of the Gauduchon metric $\eta$.

Fix a compact Hermitian manifold $(M^n, \omega)$, we consider the conformal change $\widetilde{\omega}=e^{\frac{2}{n}u}\omega$. From \cite{[Gau]}, the Chern scalar curvatures of $\widetilde{\omega}$ and $\omega$ have the following relationship:
\begin{align}\label{2.7}
-\Delta^{Ch}_\omega u+S^{Ch}(\omega)=S^{Ch}(\widetilde{\omega})e^{\frac{2}{n}u},
\end{align}
where $S^{Ch}(\omega)$ and $S^{Ch}(\widetilde{\omega})$ denote the Chern scalar curvatures of $\omega$ and $\widetilde\omega$, respectively. For convenience, let $PC(\omega)$ denote the set of $C^\infty(M)$ functions which are Chern scalar curvatures of all $\tilde\omega\in \{\omega\}$. In other words, $PC(\omega)$ is the set of $C^\infty(M)$ functions for which one can find a smooth solution of \eqref{2.7}.

At the end of this section, we recall two properties about classical Sobolev spaces on the compact Hermitian manifold $(M^n,\omega)$, which will be used later. 
\begin{lemma}\label{theorem2.3}
Assume that $n\geq2$. Then there exists a uniform constant $C$ such that 
\begin{align}
\left(\int_M|f|^{2\beta}d\mu_\omega \right)^{\frac{1}{\beta}}\leq C\left( \int_M|\nabla f|^{2}d\mu_\omega+\int_M|f|^{2}d\mu_\omega\right)
\end{align}
for any $f\in W^{1,2}(M)$, where $W^{1,2}(M)$ denotes the classical Sobolev space on $(M^n, \omega)$, and $\beta=\frac{n}{n-1}$.
\end{lemma}

\begin{lemma}[\cite{[Fon]}]\label{lemma2.4}
For any $f\in W^{k,p}(M)$ with $\int_M fd\mu_\omega=0$, $\int_M |\nabla^k f|^{p}d\mu_\omega\leq1$, $kp=2n$, there exist two constants $\delta_1=\delta_1(k,n)$ and $\delta_2=\delta_2(k, M)$ such that
\begin{align}
\int_M e^{\delta_1|f|^{\frac{p}{p-1}}} d\mu_\omega\leq \delta_2.
\end{align}
\end{lemma}

\section{The convergence of the geometric flow and prescribed non-positive Chern scalar curvatures} 
In this section, we will prove the convergence of the flow \eqref{2.1} on a compact (non-balanced) Hermitian manifold $(M^n, \omega)$ with Gauduchon degree $\Gamma(\{\omega\})<0$. 

Since the flow \eqref{2.1} is a parabolic equation with initial data $0$, it has a short-time solution $u: M\times [0,T)\rightarrow \mathbb{R}$, where $T>0$ is a constant. Firstly, we present the following maximum principle, which can be easily deduced from \cite[Lemma 4.1]{[LM]}.
\begin{lemma}\label{lemma2.2}
Assume that $u\in C^\infty(M\times[0,t_0])$ satisfies
\begin{align}
(\frac{\partial }{\partial t}-\Delta^{Ch})u\leq 0.
\end{align}
If $u(x,0)\leq c_0$, then $u(x,t)\leq c_0$ for all $(x,t)\in M\times[0,t_0]$.
\end{lemma}
In terms of the above maximum principle, we prove that the solution of \eqref{2.1} is global, i.e., $T=+\infty$.

\begin{lemma}
The flow \eqref{2.1} has a unique solution defined on $M\times [0,+\infty)$.
\end{lemma}
\proof
From $g\leq 0$, it follows that
\begin{align}\label{3.2..}
\left(\frac{\partial}{\partial t}-\Delta^{Ch}_\eta\right)\left|\frac{\partial u}{\partial t}\right|^2&=\left(\frac{\partial}{\partial t}-\Delta^{Ch}_\eta\right)\left(\Delta_\eta u-s_0+ge^{\frac{2}{n}u}\right)^2\\
&=\frac{4}{n}ge^{\frac{2}{n}u}\left(\Delta^{Ch}_\eta u-s_0+ge^{\frac{2}{n}u}\right)^2\notag\\
&\ \ \ \ -2\left|\nabla \left(\Delta^{Ch}_\eta u-s_0+ge^{\frac{2}{n}u}\right)\right|^2\notag\\
&\leq0\notag
\end{align}
for any $(x,t)\in M\times [0,T)$. On the other hand, 
\begin{align}\label{3.3..}
\left|\frac{\partial u}{\partial t}\right|^2(x,0)&=\left|\Delta^{Ch}_\eta u-s_0+ge^{\frac{2}{n}u}\right|^2(x,0)\\
&=\left|g-s_0\right|^2\leq \sup_M(|g|+|s_0|)^2<\infty.\notag
\end{align}
Using Lemma \ref{lemma2.2}, we obtain
\begin{align}
\sup_{M\times[0,T)}\left|\frac{\partial u}{\partial t}\right|\leq \sup_M(|g|+|s_0|).\label{2.5}
\end{align}
Set
\begin{align}
C_1=\sup_M(|g|+|s_0|),
\end{align}
then from \eqref{2.5}, 
\begin{align}
|u(x,t)|=|u(x,t)-u(x,0)|\leq \int^t_0\left |\frac{\partial u}{\partial t}\right|dt\leq C_1T
\end{align}
and
\begin{align}
|\Delta^{Ch}_\eta u|(x,t)\leq \sup_M|s_0|+(\sup_M|g|) e^{\frac{2}{n}C_1T}+C_1
\end{align}
for any $(x,t)\in M\times [0,T)$. According to the $L^p$-estimate for elliptic equations, we have
\begin{align}
\|u\|_{W^{2,p}(M)}\leq C(\|\Delta^{Ch}_\eta u\|_{L^p(M)}+\|u\|_{L^p(M)})\leq C_2,
\end{align}
where $C_2$ is a constant independent of $t$. By the Sobolev's embedding theorem $W^{2,p}(M)\subset C^{1}(M)$, we get $\|u\|_{C^{1}(M)}$ is bounded uniformly in $t\in [0,T)$. Hence, $\|-s_0+ge^{\frac{2}{n}u}\|_{C^\alpha(M)}$ is bounded uniformly in $t\in [0,T)$, since the Sobolev’s embedding theorem $C^1(M)\subset C^\alpha(M)$ for any $\alpha\in (0,1)$. According to \cite[Theorem 2.2.1, pp. 79-80]{[Jost]}, there exists a constant $C>0$ independent of $t$ such that
\begin{align}
\|u\|_{C^{2,\alpha}(M)}+\left\|\frac{\partial u}{\partial t}\right\|_{C^\alpha(M)}\leq C.
\end{align}
So $u(\cdot,t)\rightarrow u(\cdot, T)$ in $C^{2,\alpha}(M)$, as $t\rightarrow T$, and the local existence implies the solution $u(x,t)$ can be continued up to some time $T+\epsilon$. Therefore, the solution $u(x,t)$ is global, that is $T=+\infty$. For the uniqueness, let $u_1$ and $u_2$ be two solutions of \eqref{2.1}, then
\begin{align}
\left(\frac{\partial}{\partial t}-\Delta^{Ch}_\eta\right)(u_1-u_2)^2&=2g(u_1-u_2)\left(e^{\frac{2}{n}u_1}-e^{\frac{2}{n}u_2}\right)-2|\nabla(u_1-u_2)|^2\\
&\leq0,\notag
\end{align}
\begin{align}
(u_1-u_2)(x,0)=0.
\end{align}
According to Lemma \ref{lemma2.2}, we have $u_1\equiv u_2$.
\qed
 \begin{lemma}\label{lemma3.4}
 Let $u(x,t):M\times [0,+\infty)\rightarrow \mathbb{R}$ be the solution of the flow \eqref{2.1}. Then there exists a constant $C>0$ independent of $t$ such that 
 \begin{align}
 \|u\|_{C^0(M)}\leq C\max\{\|u\|_{L^2(M)},1\}
 \end{align}
 \end{lemma}
 \proof Multiplying the first equation of \eqref{2.1} by $u|u|^{a}$ and integrating on $M$ with respect to the Gauduchon metric $\eta$ yields 
 \begin{align}\label{3.13}
-&\int_Mu|u|^{a}\Delta^{Ch}_{\eta}ud\mu_\eta\\
&=-\int_M\left(u|u|^{a}\frac{\partial u}{\partial t}-u|u|^{a}s_0+u|u|^{a}ge^{\frac{2}{n}u}\right)d\mu_\eta\notag\\
&\leq \left(\sup_{M\times [0,\infty)}\left\vert\frac{\partial u}{\partial t}\right\vert+\sup_M|s_0|\right)\int_M |u|^{a+1}d\mu_\eta+\frac{n}{2}\int_M(-g)|u|^{a}d\mu_\eta\notag\\
&\leq C_3\left(\int_M |u|^{a+1}d\mu_\eta+\int_M|u|^{a}d\mu_\eta\right),\notag
 \end{align}
here we have used $xe^x\geq-1$ for any $x\in \mathbb{R}$, where $C_3=\max\{2C_1, \frac{n}{2}C_1\}$. On the other hand, since the map of real numbers $x\rightarrow x|x|^{\frac{a}{2}}$ is differential with derivative $(\frac{a}{2}+1)|x|^{\frac{a}{2}}$, we have
\begin{align}
\label{3.14..}
-\int_Mu|u|^{a}\Delta^{Ch}_{\eta}ud\mu_\eta&=-\int_Mu|u|^{a}\Delta_dud\mu_\eta+\int_Mu|u|^{a}(du,\theta)_\eta d\mu_\eta\\
&=\int_M\nabla(u|u|^{a})\nabla ud\mu_\eta+\frac{1}{a+2}\int_M(d(u|u|^{\frac{a}{2}})^2,\theta)_\eta d\mu_\eta\notag\\
&=\int_M(a+1)|u|^a\nabla u\cdot\nabla u d\mu_\eta+0\notag\\
&=\frac{(a+1)}{\left(\frac{a}{2}+1\right)^2}\int_M|\nabla(u|u|^{\frac{a}{2}})|^2d\mu_\eta,\notag
\end{align}
where we have used Lemma \ref{lemma2.1} and $d^*\theta=0$. Combining \eqref{3.13} and \eqref{3.14..} yields
\begin{align}\label{3.14.}
\int_M|\nabla(u|u|^{\frac{a}{2}})|^2d\mu_\eta\leq\frac{(\frac{a}{2}+1)^2}{a+1}C_3\left(\int_M |u|^{a+1}d\mu_\eta+\int_M|u|^{a}d\mu_\eta\right).
\end{align}
Using H\"older's inequality and $\text{Vol(M)}=\int_M d\mu_\eta=1$, we deduce that
\begin{align}
&\int_M|u|^{a}d\mu_\eta\leq\frac{a}{a+2} \int_M|u|^{a+2}d\mu_\eta+\frac{2}{a+2},\label{3.16...}\\
&\int_M|u|^{a+1}d\mu_\eta\leq\frac{a+1}{a+2} \int_M|u|^{a+2}d\mu_\eta+\frac{1}{a+2}.\label{3.17...}
\end{align}
Hence, 
\begin{align}
\int_M|\nabla(u|u|^{\frac{a}{2}})|^2d\mu_\eta\leq (a+1)C_3\int_M|u|^{a+2}d\mu_\eta+\frac{3}{2}C_3.\label{3.17.}
\end{align}
Now applying the Sobolev's inequality (cf. Lemma \ref{theorem2.3}) to $f=u|u|^{\frac{a}{2}}$, we obtain
\begin{align}
\left(\int_M|u|^{(a+2)\beta}d\mu_\eta\right)^{\frac{1}{\beta}}\leq C_4\left(\int_M|\nabla(u|u|^{\frac{a}{2}})|^2d\mu_\eta+\int_M|u|^{a+2}d\mu_\eta \right),\label{3.15.}
\end{align}
where $C_4$ is a constant independent of $u$, and $\beta=\frac{n}{n-1}$. From \eqref{3.17.} and \eqref{3.15.}, it follows that
\begin{align}
\left(\int_M|u|^{(a+2)\beta}d\mu_\eta\right)^{\frac{1}{\beta}}&\leq C'_3C_4(a+2)\int_M|u|^{a+2}d\mu_\eta+\frac{3}{2}C_3C_4,\label{3.19}\\
&\leq 2C'_3C_4 (a+2)\left(\int_M|u|^{a+2}d\mu_\eta+ 1\right)\notag\\
&\leq C_5(a+2)\max\left\{\int_M|u|^{a+2}d\mu_\eta,1\right\}\notag
\end{align}
where $C'_3=\max\{C_3,1\}$, $C_5=4C'_3C_4$. Set $p=a+2\geq2$, then from \eqref{3.19} we have
\begin{align}
\max\{\|u\|_{L^{p\beta}(M)},1\}\leq (C'_5p)^{\frac{1}{p}}\max\{\|u\|_{L^p(M)},1\},
\end{align}
where $C'_5=\max\{C_5,1\}$. Therefore, for any positive integer $k$, 
\begin{align}
\max\{\|u\|_{L^{p\beta^{k+1}}(M)},1\}&\leq (C'_5)^{\frac{1}{p\beta^k}}(p\beta^k)^{\frac{1}{p\beta^k}}\max\{\|u\|_{L^{p\beta^k}(M)},1\}\\
&\leq\dots\notag\\
&\leq (C'_5p)^{\frac{1}{p}\sum_{i=0}^k\frac{1}{\beta^i}}\beta^{\frac{1}{p}\sum_{i=0}^k\frac{i}{\beta^i}}\max\{\|u\|_{L^p(M)},1\}.\notag
\end{align}
Since $\beta=\frac{n}{n-1}>1$, we see that $\sum_{i=0}^\infty\frac{1}{\beta^i}<+\infty$ and $\sum_{i=0}^\infty\frac{i}{\beta^i}<+\infty$. Then there exists a constant $C_6$ independent of $t$ and $k$ such that
\begin{align}
\max\{\|u\|_{L^{p\beta^{k+1}}(M)},1\}\leq C_6\max\{\|u\|_{L^p(M)},1\}.
\end{align}
Pick $p=2$ and let $k\rightarrow +\infty$ in the above inequality, we obtain
\begin{align}
\|u\|_{C^0(M)}\leq\max\{\|u\|_{C^0(M)},1\}\leq C_6\max\{\|u\|_{L^2(M)},1\}.\label{3.23.}
\end{align}
\qed
\begin{lemma}
If $(M^n,\omega)$ is of negative Gauduchon degree, then there exists a constant $C_7>0$ independent of $t$ such that
\begin{align}
\|u\|_{L^2(M)}<C_7,
\end{align} 
where $u$ is the solution of \eqref{2.1}. Therefore, from Lemma \ref{lemma3.4}, we have 
\begin{align}\label{3.26...}
\|u\|_{C^0(M)}\leq C'_7. 
\end{align}
where $C'_7=C_6\max\{C_7,1\}$.
\end{lemma}
\proof We will prove it by contradiction. Assume that $\|u\|_{L^2(M)}$ is not bounded uniformly with respect to $t\in[0,+\infty)$, then there is a sequence $\{t_i\}_{i=1}^\infty$ with $\lim_{i\rightarrow +\infty}t_i=+\infty$ such that 
\begin{align}
\lim_{i\rightarrow +\infty}\|u(\cdot,t_i)\|_{L^2(M)}=+\infty,
\end{align}
and 
\begin{align}
\frac{d}{dt}\vert_{t=t_i}\|u(\cdot,t)\|_{L^2(M)}\geq 0\label{3.26.}
\end{align}
for $t_i$ large enough.
Set 
\begin{align}
u_i=u(\cdot,t_i),\ l_i=\|u(\cdot,t_i)\|_{L^2(M)},\ w_i=\frac{u_i}{l_i},
\end{align}
then $\|w_i\|_{L^2(M)}=1$. Moreover, according to \eqref{3.23.}, we obtain
\begin{align}
\|w_i\|_{C^0(M)}\leq C_8,
\end{align}
where $C_8$ is a positive constant independent of $i$. Using $d^*\theta=0$, \eqref{2.1} and H\"older's inequality, we have the following computation
\begin{align}
\int_M|\nabla u_i|^2d\mu_\eta&=-\int_Mu_i\Delta^{Ch}_\eta u_id\mu_\eta\label{3.28}\\
&=\int_M u_i\left(-s_0+ge^{\frac{2}{n}u_i}-\frac{\partial u}{\partial t}\vert_{t=t_i}\right)d\mu_\eta\notag\\
&\leq \left(\sup_M |s_0|+\sup_{M\times[0,+\infty)}\left|\frac{\partial u}{\partial t}\right|\right)\int_M |u_i|d\mu_\eta+\frac{n}{2}\int_M(-g)d\mu_\eta\notag\\
&\leq C_{9}\|u_i\|_{L^2(M)}+\frac{n}{2}\int_M(-g)d\mu_\eta,\notag
\end{align}
where $C_{9}=2C_1$, and we have used  $xe^x\geq-1$ for any $x\in \mathbb{R}$. Dividing both sides of \eqref{3.28} by $l_i^2$ yields
\begin{align}
\int_M|\nabla w_i|^2d\mu_\eta\leq \frac{C_{9}}{l_i}+\frac{n}{2}\frac{\int_M(-g)d\mu_\eta}{l_i^2},
\end{align}
so
\begin{align}
\lim_{\ i\rightarrow +\infty}\int_M|\nabla w_i|^2d\mu_\eta=0.\label{3.30}
\end{align}
Hence, $\{w_i\}_{i=1}^\infty$ is bounded in $W^{1,2}(M)$. Combining the Sobolve's embedding theorem $W^{1,2}(M)\subset\subset L^2(M)$, there exists a subsequence of $\{w_i\}_{i=1}^\infty$, also denoted by $\{w_i\}_{i=1}^\infty$, and $w_{\infty}\in W^{1,2}(M)$ such that
\begin{align}
&w_i\rightarrow w_{\infty}\ \text{in\ } L^2(M),\\
&w_i\rightharpoonup w_{\infty}\ \text{weakly\ in\ } W^{1,2}(M).\label{3.32}
\end{align}
From \eqref{3.30}-\eqref{3.32} and $\|w_i\|_{L^2(M)}=1$, it follows that
\begin{align}
w_\infty=C^*\not=0\ \text{almost\ everywhere\ in\ } M,
\end{align} 
where $C^*$ is a constant. By \eqref{2.1}, we have
\begin{align}
\int_M(-g)e^{\frac{2}{n}u_i}d\mu_\eta&=\int_M\Delta^{Ch}_\eta u_id\mu_\eta-\int_Ms_0d\mu_\eta-\int_M\frac{\partial u}{\partial t}|_{t=t_i}d\mu_\eta\\
&\leq -\int_Ms_0d\mu_\eta+C_1,\notag
\end{align}
thus,
\begin{align}\label{3.35}
\frac{2}{n}\int_M(-g)w_id\mu_\eta&\leq\int_M(-g)\frac{e^{\frac{2}{n}l_iw_i}}{l_i}d\mu_\eta\\
&\leq\frac{-\int_Ms_0d\mu_\eta+C_1}{l_i}\rightarrow 0, \text{as}\ i\rightarrow +\infty.\notag
\end{align}
On the other hand, according to \eqref{3.32}, we obtain
\begin{align}
\lim_{i\rightarrow +\infty}\int_M(-g)w_id\mu_\eta=C^*\int_M(-g)d\mu_\eta.\label{3.36}
\end{align}
Combining \eqref{3.35} and \eqref{3.36} yields $C^*<0$, because of $g\leq 0\ (\not\equiv 0)$ and $C^*\not=0$. Now in terms of \eqref{2.1}, \eqref{3.32}, we perform the following computation
\begin{align}
\frac{d}{dt}|_{t=t_i}\|u\|_{L^2(M)}&=\frac{1}{l_i}\int_Mu_i\frac{\partial u}{dt}|_{t=t_i}d\mu_\eta\\
&=\frac{1}{l_i}\int_Mu_i\left(\Delta^{Ch}_\eta u_i-s_0+ge^{\frac{2}{n}u_i} \right)d\mu_\eta\notag\\
&\leq-\frac{1}{l_i}\int_M|\nabla u_i|^2d\mu_\eta-\int_Mw_is_0d\mu_\eta+\frac{n\int_M(-g)d\mu_\eta}{2l_i}\notag\\
&\leq-\int_Mw_is_0d\mu_\eta+\frac{n\int_M(-g)d\mu_\eta}{2l_i}\rightarrow -C^*\int_Ms_0d\mu_\eta,\notag
\end{align}
as $i\rightarrow +\infty$, which is a contradiction with \eqref{3.26.}, since Gauduchon degree $\Gamma(\{\omega\})=\int_Ms_0d\mu_\eta<0$ and $C^*<0$.
\qed

In terms of the above lemmas, we can prove the convergence of the flow \eqref{2.1} now.
\begin{theorem}
Let $(M^n,\omega)$ be a compact Hermitian manifold with complex dimension $n\geq 2$ and Gauduchon degree $\Gamma (\{\omega\})<0$, and assume that $g\in C^\infty(M)$. Suppose either
\begin{enumerate}
\item $g$ is negative on $M$
\end{enumerate}
or
\begin{enumerate}[\ \ \ \ (2)]
\item $g$ is nonpositive and not equal to zero identically on $M$ and $\{\omega\}$ contains a balanced metric,
\end{enumerate}
then there is a subsequence $t_i\rightarrow +\infty$ such that $u(\cdot,t_i)$ converges to $u_\infty$ in $C^{1,\alpha}(M)$, where $\alpha\in (0,1)$, and $u_\infty$ is the unique smooth solution of \eqref{1.6.}.
\end{theorem}

\proof (1) Since $g<0$ on $M$, then by \eqref{3.2..} we have
\begin{align}\label{3.41}
\left(\frac{\partial}{\partial t}-\Delta^{Ch}_\eta\right)\left|\frac{\partial u}{\partial t}\right|^2\leq -C_{10}\left|\frac{\partial u}{\partial t}\right|^2,
\end{align}
where $C_{10}=-\frac{4}{n}\left(\max_M g\right)e^{\frac{2}{n}C'_7}>0$. Set
\begin{align}
\varphi(x,t)=e^{C_{10}t}\left|\frac{\partial u}{\partial t}\right|^2.
\end{align}
From \eqref{3.41} and \eqref{3.3..}, it follows that
\begin{equation}
\left\{
\begin{aligned}
& \left(\frac{\partial}{\partial t}-\Delta^{Ch}_\eta\right)\varphi\leq 0,\\
&\varphi(x,0)=|-s_0+g|^2.
\end{aligned}
\right.
\end{equation}
Using Lemma \ref{lemma2.2} yields
\begin{align}
\varphi(x,t)=e^{C_{10}t}\left|\frac{\partial u}{\partial t}\right|^2\leq \max_M|-s_0+g|^2,
\end{align}
so 
\begin{align}\label{3.45}
\max_{x\in M}\left|\frac{\partial u}{\partial t}\right|^2(x,t)\leq e^{-C_{10}t}\max_M|-s_0+g|^2\rightarrow 0,\ \text{as}\ t\rightarrow +\infty.
\end{align}
On the other hand, According to \eqref{2.1}, \eqref{2.5} and \eqref{3.26...}, we obtian
\begin{align}
\|\Delta^{Ch}_\eta u\|_{C^0(M)}\leq \sup_M|s_0|+(\sup_M|g|) e^{\frac{2}{n}C'_7}+C_1,
\end{align}
and thus $\|u\|_{W^{2,p}(M)}$ is bounded uniformly with respect to $t\in [0, +\infty)$ for any $p>1$, because of the $L^p$-estimate for elliptic equations. Therefore, by the Sobolev embedding theorem $W^{2,p}(M)\subset\subset C^{1,\alpha}(M)$ for any $p>n$ and $\alpha\in (0,1)$, there exists a subsequence $t_i\rightarrow +\infty$ such that
\begin{align}\label{3.47}
u(\cdot,t_i)\rightarrow u_\infty, \text{\ in\ } C^{1,\alpha}(M).
\end{align}
Clearly, \eqref{3.45} and \eqref{3.47} implies $u_\infty$ is a weak solution of \eqref{1.3}. By the regularity results for elliptic equations and the maximum principle, we can conclude that $u_\infty$ is the unique smooth solution of \eqref{1.3}.

(2) Since $\{\omega\}$ contains a balanced metric, then from Theorem \ref{theorem2.2} it follows that the Gauduchon metric $\eta$ in \eqref{2.1} is balanced. By a similar computation, we deduce that
\begin{align}
\int_M\left|\frac{\partial u}{\partial t}\right|^2d\mu_\eta&=\int_M\left(\Delta^{Ch} u-s_0+ge^{\frac{2}{n}u}\right)^2d\mu_\eta\\
&=\frac{d}{dt}\int_M\left(us_0+\frac{1}{2}|\nabla u|^2-\frac{n}{2}\left(e^{\frac{2}{n}u}-1\right)g\right)d\mu_\eta.\notag
\end{align}
Since $u(x,0)=0$, we have
\begin{align}\label{3.49}
\int_0^s\int_M\left|\frac{\partial u}{\partial t}\right|^2d\mu_\eta dt=\int_M\left(us_0+\frac{1}{2}|\nabla u|^2-\frac{n}{2}\left(e^{\frac{2}{n}u}-1\right)g\right)d\mu_\eta
\end{align}
for any $s\geq 0$. According to the proof of (1) in Theorem \ref{theorem1.1}, we have proved that $\|u\|_{C^{1,\alpha}(M)}$ is uniformly bounded with respect to $t\in [0,+\infty)$. Hence, by \eqref{3.49}, there exists a constant $C>0$ such that
\begin{align}
\int_0^{+\infty}\int_M\left|\frac{\partial u}{\partial t}\right|^2d\mu_\eta dt<C<+\infty,
\end{align}
which implies that there exists a subsequence $t_i\rightarrow +\infty$ such that
\begin{align}
\int_M\left|\frac{\partial u}{\partial t}\right|^2(\cdot, t_i)d\mu_\eta\rightarrow 0,
\end{align}
as $t_i\rightarrow +\infty$. Then according to the same argument as the proof of (1) in Theorem \ref{theorem1.1}, we may also obtain a subsequence $t_i\rightarrow +\infty$ such that $u(\cdot,t_i)\rightarrow u_\infty$ in $C^{1,\alpha}(M)$, where $u_\infty$ is the unique smooth solution of \eqref{1.3}.
\qed

\begin{remark}
In \cite{[WZ]}, Y. Wang and X. Zhang established an existence theorem of a Kazdan-Warner type equation on a class of non-compact Riemannian manifolds by using the heat flow method. Later, R. Wang \cite{[Wan]} solved the Dirichlet problem of a Kazdan-Warner type equation on compact Hermitian manifolds with smooth boundary by combining the heat flow method and the solvability of Poisson's equations, and extended Wang-Zhang's existence theorem to a similar class of non-compact Gauduchon manifolds.
\end{remark}

\section{Prescribed sign-changing Chern scalar curvatures}
In this section, we will consider the case that the candidate curvature function $g$ is sign-changing. For any sign-changing function $g\in C^\infty(M)$ on the compact Hermitian manifold $M$, there always exist a non-constant smooth function $g_0$ satisfying $\max_M g_0=0$ and a constant $\lambda>0$ such that $g=g_0+\lambda$. In terms of the implicit function theorem, we obtain the following existence and non-existence results.

\begin{theorem}\label{theorem3.8}
Let $(M^n,\omega)$ be a compact Hermitian manifold with complex dimension $n\geq 2$ and Gauduchon degree $\Gamma (\{\omega\})<0$. Let $g_0$ be a non-constant smooth function on $M$ with $\max_M g_0=0$. Then there exists a constant $\lambda^*\in (0,-\min_M g_0]$ such that 
\begin{enumerate}
\item If $0<\lambda<\lambda^*$, then $g_0+\lambda\in PC(\omega)$.
\item If $\lambda>\lambda^*$, then $g_0+\lambda\not\in PC(\omega)$.
\end{enumerate}
\end{theorem}
\proof Consider the Banach manifolds
\begin{align}
&X=\left\{(u,\lambda)\in C^\infty(M)\times [0,+\infty) : \int_M s_0d\mu_\eta=\int_M (g_0+\lambda)e^{\frac{2}{n}u}d\mu_\eta \right\},\\
&Y=\left\{w\in C^\infty(M): \int_M wd\mu_\eta=0\right\},
\end{align}
and the map $F:X\rightarrow Y$ given by
\begin{align}
F(u,\lambda)=-\Delta_\eta^{Ch}u+s_0-(g_0+\lambda)e^{\frac{2}{n}u},
\end{align}
where $\eta$ is the unique Gauduchon representative in the conformal class $\{\omega\}$ with volume $1$, and $s_0$ denotes the corresponding Chern scalar curvature. According to \cite[Theorem 2.5]{[Fus]}, there exists a unique solution $u_0\in C^\infty(M)$ of $-\Delta_\eta^{Ch}u+s_0-g_0e^{\frac{2}{n}u}=0$, i.e., $F(u_0,0)=0$. Moreover, we have 
\begin{align}
\frac{\partial F}{\partial u}|_{(u_0,0)}v=\frac{d}{dt}|_{t=0}F(u_0+tv,0)=-\Delta_\eta^{Ch}v-\frac{2}{n}g_0e^{\frac{2}{n}u_0}v,
\end{align}
which is invertible on $C^\infty(M)$. Indeed, if $-\Delta_\eta^{Ch}v-\frac{2}{n}g_0e^{\frac{2}{n}u_0}v=0$, then
\begin{align}
0\geq\int_M\frac{2}{n}g_0e^{\frac{2}{n}u_0}v^2d\mu_\eta&=-\int_Mv\Delta_\eta^{Ch}vd\mu_\eta\\
&=-\int_Mv\Delta_dvd\mu_\eta+\frac{1}{2}\int_M(dv^2,\theta)d\mu_\eta\notag\\
&=\int_M|\nabla v|^2d\mu_\eta,\notag
\end{align}
where the last equality holds since $\eta$ is a Gauduchon metric, i.e., $d^*\theta=0$. Hence, $v$ is constant and thus $v\equiv 0$. Applying the implicit function theorem, there exists a constant $\epsilon_0>0$ and smooth map $(u,id) : [0, \epsilon_0)\times [0,\epsilon_0)\rightarrow X$ such that
\begin{align}
F(u(\lambda),\lambda)=0,\ \forall \ \lambda\in [0, \epsilon_0),
\end{align}
that is
\begin{align}
-\Delta_\eta^{Ch}u(\lambda)+s_0=(g_0+\lambda)e^{\frac{2}{n}u(\lambda)},\ \forall \ \lambda\in [0, \epsilon_0).
\end{align}
Set 
\begin{align}
\lambda^*=\sup_M\{\lambda>0: -\Delta_\eta^{Ch}u+s_0=(g_0+\lambda)e^{\frac{2}{n}u}\ \text{is\ solvable}\}>0.
\end{align}
Note that $0<\lambda^*\leq-\min_M g_0<+\infty$. Indeed, If $-\Delta_\eta^{Ch}u+s_0=(g_0+\lambda)e^{\frac{2}{n}u}$ is solvable, then 
\begin{align}
\int_M (g_0+\lambda)e^{\frac{2}{n}u} d\mu_\eta=\int_M s_0d\mu_\eta=\Gamma(\{w\})<0,
\end{align}
which implies $g_0+\lambda$ is negative somewhere on $M$, so $\lambda<-\min_M g_0$, therefore $\lambda^*\leq -\min_M g_0$. If $0<\lambda<\lambda^*$, then there is a constant $\lambda_0\in (\lambda, \lambda^*)$ such that $-\Delta_\eta^{Ch}u+s_0=(g_0+\lambda_0)e^{\frac{2}{n}u}$ is solvable. From \cite[Proposition 2.6 (b)]{[Fus]}, it follows that $g_0+\lambda\in PC(\omega)$, which completes the proof of $(1)$ in Theorem \ref{theorem3.8}. For the assertion $(2)$ in Theorem \ref{theorem3.8}, it is obvious due to the definition of $\lambda^*$.
\qed

Now we consider the case that $\lambda=\lambda^*$. For arbitrary Hermitian manifold $(M^n,\omega)$, whether $g_0+\lambda^*$ belongs to $PC(\omega)$ is still open. However, when $\dim_{\mathbb{C}}M=2$, we prove the following theorem.

\begin{theorem}\label{theorem4.2}
Let $(M^n,\omega)$ be a compact Hermitian manifold with complex dimension $n=2$ and Gauduchon degree $\Gamma (\{\omega\})<0$, and suppose that $\{\omega\}$ contains a balanced metric $\eta$ with volume $1$. Let $g_0$ be a non-constant smooth function on $M$ with $\max_M g_0=0$, and $\lambda^*$ be a constant as in Theorem \ref{theorem3.8}. Then $g_0+\lambda^*\in PC(\omega)$. 
\end{theorem}

According to Theorem \ref{theorem3.8}, for any $\lambda\in (0,\lambda^*)$, there always exits a solution $u_\lambda\in C^\infty(M)$ of
\begin{align}
-\Delta_\eta^{Ch}u_\lambda+s_0=(g_0+\lambda)e^{\frac{2}{n}u_\lambda}.\label{3.49.}
\end{align}
Note that the solution $u_\lambda$ is usually not unique. In order to prove that \eqref{3.49.} also has a solution for $\lambda=\lambda^*$, we will divide the proof into the following several steps. Remark that here we combine the sprit of Ding-Liu's argument \cite{[DL]} and the Moser iteration method to prove it. 
\begin{lemma}\label{lemma4.3.}
There exists a function $u_\lambda\in C^\infty(M)\cap X$ such that 
\begin{align}
I_\lambda(u_\lambda)=\inf_{u\in X}I_\lambda(u)=\inf_{u\in X}\int_M \left(|\nabla u|^2+2s_0u-n(g_0+\lambda)e^{\frac{2}{n}u}\right)d\mu_\eta,
\end{align}
where $X=\{u\in W^{1,2}(M): u_1\leq u\leq u_2\ \text{a.e.\ on}\ M\}$, $u_1$ (Resp. $u_2$) is a smooth sub-solution (Resp. super-solution) of \eqref{3.49.} with $u_1<u_2$. (Note that both $u_1$ and $u_2$ are related to $\lambda$.) Therefore,
\begin{align}
&0=\frac{d}{dt}|_{t=0}I_\lambda(u_\lambda+t\varphi)=2\int_M\left(\nabla u_\lambda\cdot \nabla\varphi+s_0\varphi-(g_0+\lambda)e^{\frac{2}{n}u_\lambda}\varphi\right)d\mu_\eta,\label{4.12...}\\
&0\leq\frac{d^2}{dt^2}|_{t=0}I_\lambda(u_\lambda+t\varphi)=2\int_M|\nabla \varphi|^2d\mu_\eta-\frac{4}{n}\int_M (g_0+\lambda)e^{\frac{2}{n}u_\lambda}\varphi^2 d\mu_\eta\label{3.52.}
\end{align}
for any $\varphi\in W^{1,2}(M)$.
\end{lemma}
\begin{remark}
$u_1, u_2$ may be chosen similar to \cite{[DL]}: let $u_2$ be a smooth solution of $-\Delta^{Ch}_\eta u+s_0=(g_0+\lambda_1)e^{\frac{2}{n}u}$, where $\lambda_1\in (\lambda, \lambda^*)$, then it is easy to see that $u_2$ is a super-solution of \eqref{3.49.}. Since $\Gamma(\{\omega\})<0$, there is a sub-solution $u_1$ of \eqref{3.49.} with $u_1<u_2$ (see the formula \eqref{4.16} below). 
\end{remark}
\proof The existence of $u_\lambda$ is due to the standard variational theory, and the smoothness of $u_\lambda$ follows from the regularity results for elliptic equations. Indeed, it is easy to see that $X$ is a convex, (weakly) closed subset of the Hilbert space $W^{1,2}(M)$, and $I_\lambda: X\rightarrow \mathbb{R}\cup +\infty$ is coercive and sequentially weakly lower semi-continuous on $X$ (because of Tonelli-Morrey Theorem, see e.g. \cite[Theorem 1.6, p. 9]{[Str]}). 
\qed
\begin{lemma}\label{lemma4.4}
There exists a constant $C=C(\lambda^*, s_0, g_0,\Gamma(\{\omega\}), M)>0$ such that
\begin{align}
u_\lambda\geq-C.\label{3.53}
\end{align}
\end{lemma}
\proof This proof is similar to that in \cite{[DL]}, for the convenience of readers, we still give the details. Let $\phi_t=f-t$, where $t$ is a constant and $f$ satisfies 
\begin{align}
\Delta^{Ch}_\eta f=s_0-\int_M s_0d\mu_\eta.
\end{align}
Then for any $t\geq t_0=-\frac{n}{2}\log{\frac{\Gamma(\{\omega\})e^{\frac{2}{n}\|f\|_{C^0(M)}}}{\|g_0\|_{C^0(M)}+\lambda^*}}+1$ we have
\begin{align}\label{4.16}
-&\Delta^{Ch}_\eta\phi_t+s_0-(g_0+\lambda)e^{\frac{2}{n}\phi_t}\\
&\leq\Gamma(\{\omega\})+(\|g_0\|_{C^0(M)}+\lambda^*)e^{\frac{2}{n}\left(\|f\|_{C^0(M)}-t_0\right)}< 0.\notag
\end{align}
Now we claim that $u_\lambda\geq \phi_{t_0}$ on $M$. Otherwise, there is a constant $t_1\in (t_0, +\infty)$ such that 
\begin{align}
u_\lambda\geq \phi_{t_1}\ \text{on}\ M, \ \ \ u_\lambda(x_0)=\phi_{t_1}(x_0)\ \text{for some}\ x_0\in M.
\end{align}
According to the maximum principle, we get $\Delta^{Ch}_\eta(u_\lambda-\phi_{t_1})(x_0)\geq 0$. Then at $x_0$
\begin{align}
-&\Delta^{Ch}_\eta\phi_{t_1}+s_0-(g_0+\lambda)e^{\frac{2}{n}\phi_{t_1}}\\
&=\Delta^{Ch}_\eta(u_\lambda-\phi_{t_1})+\left(-\Delta^{Ch}_\eta u_\lambda+s_0-(g_0+\lambda)e^{\frac{2}{n}u_\lambda}\right)\notag\\
&\ \ \ \ \ \ \ +(g_0+\lambda)\left(e^{\frac{2}{n}u_\lambda}-e^{\frac{2}{n}\phi_{t_1}}\right)\notag\\
&=\Delta^{Ch}_\eta(u_\lambda-\phi_{t_1})\geq 0,\notag
\end{align}
which is a contradiction with \eqref{4.16}. From $u_\lambda\geq \phi_{t_0}$ on $M$, \eqref{3.53} follows.
\qed

\begin{lemma}
$M_-=\{x\in M: g_0+\lambda^*<0\}$ is not empty.
\end{lemma}
\proof
Set $v_\lambda=u_\lambda-v_0$, where $v_0$ is the solution of $-\Delta_\eta^{Ch} v_0+s_0=-e^{\frac{2}{n}v_0}$. Then we have
\begin{align}
-\Delta_\eta^{Ch} v_\lambda=(g_0+\lambda)e^{\frac{2}{n}u_\lambda}+e^{\frac{2}{n}v_0}\geq (g_0+\lambda)e^{\frac{2}{n}u_\lambda},
\end{align}
thus,
\begin{align}
\int_M(g_0+\lambda^*)e^{\frac{2}{n}v_0}d\mu_\eta&=\lim_{\lambda\rightarrow\lambda^*_-}\int_M(g_0+\lambda)e^{\frac{2}{n}v_0}d\mu_\eta\\
&\leq-\int_Me^{-\frac{2}{n}v_\lambda}\Delta_\eta^{Ch} v_\lambda d\mu_\eta\notag\\
&=-\frac{2}{n}\int_Me^{-\frac{2}{n}v_\lambda}|\nabla v_\lambda|^2d\mu_\eta\leq 0.\notag
\end{align}
Hence, $g_0+\lambda^*$ is negative somewhere in $M$, i.e., $M_-$ is not empty.
\qed

\begin{lemma}\label{lemma4.6}
$u_\lambda^+(x)=\max\{u_\lambda(x),0\}$ is locally uniformly $W^{1,2}$-bounded in the open set $M_-$ as $\lambda\rightarrow \lambda^*$.
\end{lemma}
\proof To see this, it is sufficient to prove for any $p\in M_-$, there exists a constant $C>0$ independent of $\lambda$ such that
\begin{align}
\|u^+_\lambda\|_{W^{1,2}\left(B_{\frac{d}{2^{n_0+1}}}\right)}\leq C,\label{3.55}
\end{align}
where $d=dist_\eta(p,\partial M_-)$, $n_0$ is an integer depending only on $p$ and is large enough so that the closed geodesic ball $\overline{B_{\frac{d}{2^{n_0+1}}}}=\{x\in M: dist_\eta(x,p)\leq\frac{d}{2^{n_0+1}} \}$ is contained within the cut locus of $p$. In order to show \eqref{3.55}, we choose a smooth cut-off function $\xi$ on $M$ satisfying
\begin{align}\label{4.22...}
\xi(x)=
\begin{cases}
1& x\in B_{\frac{d}{2^{n_0+1}}}\\
0& x\in M\setminus B_{\frac{d}{2^{n_0}}} 
\end{cases}
,\ 0\leq\xi\leq1,\ |\nabla \xi|\leq \frac{A}{d}\xi^{\frac{1}{2}},
\end{align}
where $A>0$ is a constant related to $n_0$. (Note that the cut-off function $\xi$ is different from the one used in \cite{[DL]}.) Indeed, we first take a cut-off function $\phi\in C^\infty([0,+\infty))$ with
\begin{align}
\phi|_{[0,1]}=1,\ \phi|_{[2,+\infty)}=0,\ -A'|\phi|^{\frac{1}{2}}\leq \phi'\leq0,
\end{align}
where $A'>0$ is a constant. Set $\xi(x)=\phi(\frac{2^{n_0+1}}{d}r(x))\ (\forall x\in M)$, where $r(x)$ is the Riemannian distance between $p$ and $x$ with respect to $\eta$. Then it is easy to see that the above $\xi$ satisfies \eqref{4.22...}. Pick $\varphi=\xi^2u_\lambda^+$ in \eqref{4.12...}, we get
\begin{align}
\int_{B_{\frac{d}{2^{n_0}}}}\left(\nabla u^+_\lambda\cdot\nabla (\xi^2u^+_\lambda)+s_0\xi^2u^+_\lambda-(g_0+\lambda)e^{\frac{2}{n}u^+_\lambda}\xi^2u^+_\lambda\right)d\mu_\eta=0.\label{3.57}
\end{align}
Making use of
\begin{align}
&\nabla u^+_\lambda\cdot\nabla (\xi^2u^+_\lambda)=|\nabla (\xi u^+_\lambda)|^2-|\nabla \xi|^2(u^+_\lambda)^2,\label{4.25...}\\
&g_0+\lambda\leq g_0+\lambda^*\leq-\epsilon<0  \ \text{for some} \ \epsilon=\epsilon(\lambda^*,d)\in (0,1),
\end{align}
(compare \eqref{4.25...} with \cite[page 1062, line -9]{[DL]}), it follows from \eqref{3.57} that
\begin{align}
0\geq\int_{B_{\frac{d}{2^{n_0}}}}\left(|\nabla (\xi u^+_\lambda)|^2-|\nabla \xi|^2(u^+_\lambda)^2+s_0\xi^2u^+_\lambda+\epsilon e^{\frac{2}{n}u^+_\lambda}\xi^2u^+_\lambda\right)d\mu_\eta.\label{3.60}
\end{align}
In terms of \eqref{4.22...} and Cauchy's inequality, we deduce that
\begin{align}\label{3.61}
\frac{d^2}{A^2} \int_{B_{\frac{d}{2^{n_0}}}}|\nabla \xi|^2(u^+_\lambda)^2d\mu_\eta&\leq\int_{B_{\frac{d}{2^{n_0}}}}\xi(u^+_\lambda)^2d\mu_\eta\\
&\leq\delta_1\int_{B_{\frac{d}{2^{n_0}}}}\xi^2(u^+_\lambda)^4d\mu_\eta+(4\delta_1)^{-1}\text{Vol}({B_{\frac{d}{2}}})\notag
\end{align}
and
\begin{align}
\left|\int_{B_{\frac{d}{2^{n_0}}}}s_0\xi^2u^+_\lambda d\mu_\eta\right|\leq (4\delta_2)^{-1}\int_{B_{\frac{d}{2^{n_0}}}}s_0^2d\mu_\eta+\delta_2\int_{B_{\frac{d}{2^{n_0}}}}\xi^2(u^+_\lambda)^2d\mu_\eta\label{3.62}
\end{align}
for any $\delta_1,\delta_2>0$. From \eqref{3.60}-\eqref{3.62} and the facts $e^{2t}\geq t^3$, $t^4>t^2-1\ (\forall\ t\in \mathbb{R})$, we obtain
\begin{align}
0& \geq \int_{B_{\frac{d}{2^{n_0}}}}|\nabla (\xi u^+_\lambda)|^2d\mu_\eta+\left(\frac{\epsilon}{n^3}-\frac{A^2}{d^2}\delta_1\right)\int_{B_{\frac{d}{2^{n_0}}}}\xi^2(u^+_\lambda)^4d\mu_\eta\label{3.63}\\
&\ \ \ \ \ \ -\frac{A^2}{4\delta_1d^2}\text{Vol}({B_{\frac{d}{2^{n_0}}}})-(4\delta_2)^{-1}\int_{B_{\frac{d}{2^{n_0}}}}s_0^2d\mu_\eta-\delta_2\int_{B_{\frac{d}{2^{n_0}}}}\xi^2(u^+_\lambda)^2d\mu_\eta\notag\\
&\geq \int_{B_{\frac{d}{2^{n_0}}}}|\nabla (\xi u^+_\lambda)|^2d\mu_\eta+\left(\frac{\epsilon}{n^3}-\frac{A^2}{d^2}\delta_1-\delta_2\right)\int_{B_{\frac{d}{2^{n_0}}}}\xi^2(u^+_\lambda)^2d\mu_\eta\notag\\
&\ \ \ \ \ \ -\left(\frac{1}{4\delta_1}\frac{A^2}{d^2}+\frac{\epsilon}{n^3}-\frac{A^2}{d^2}\delta_1\right)\text{Vol}({B_{\frac{d}{2^{n_0}}}})-(4\delta_2)^{-1}\int_{B_{\frac{d}{2^{n_0}}}}s_0^2d\mu_\eta.\notag
\end{align}
Picking $\delta_1=\frac{\epsilon d^2}{4A^2n^3}, \delta_2=\frac{\epsilon}{4n^3}$ in \eqref{3.63} yields
\begin{align}
\int_{B_{\frac{d}{2^{n_0}}}}|&\nabla (\xi u^+_\lambda)|^2d\mu_\eta+\frac{\epsilon}{2n^3}\int_{B_{\frac{d}{2^{n_0}}}}\xi^2(u^+_\lambda)^2d\mu_\eta\\
&\leq\frac{n^3}{\epsilon}\int_{B_{\frac{d}{2^{n_0}}}}s_0^2d\mu_\eta+\left(\frac{n^3}{\epsilon}\frac{A^4}{d^4}+\frac{3\epsilon}{4n^3}\right)\text{Vol}({B_{\frac{d}{2^{n_0}}}}),\notag
\end{align}
which implies that
\begin{align}
\| u^+_\lambda\|^2_{W^{1,2}\left(B_{\frac{d}{2^{n_0+1}}}\right)}\leq \frac{2n^6}{\epsilon}\int_{B_{\frac{d}{2^{n_0}}}}s_0^2d\mu_\eta+(\frac{2n^6}{\epsilon^2}\frac{A^4}{d^4}+\frac{3}{2})\text{Vol}({B_{\frac{d}{2^{n_0}}}}).
\end{align}
\qed

\begin{lemma}\label{lemma4.7}
$u_\lambda$ is locally uniformly $W^{1,2}$-bounded in the open set $M_-$ as $\lambda\rightarrow \lambda^*$.
\end{lemma}
\proof Let $D$ be an arbitrary open subset in $M_-$ with $D\subset\subset M_-$. According to \eqref{3.53} and Lemma \ref{lemma4.6}, we have the following
\begin{align}\label{3.67..}
\int_D |u_\lambda|^2d\mu_\eta&=\int_{D\cap\{x\in D : -C\leq u_\lambda<0\}} |u_\lambda|^2d\mu_\eta+\int_{D\cap\{x\in D : u_\lambda\geq0\}} (u_\lambda^+)^2d\mu_\eta\\
&\leq C^2\text{Vol(D)}+C^*,\notag
\end{align}
where $C$ and $C^*$ are two constants independent of $\lambda$. Furthermore, if we take the test function $\varphi=\xi^2u_\lambda$ in \eqref{4.12...}, then 
\begin{align}
\int_{M}\left(\nabla u_\lambda\cdot\nabla (\xi^2u_\lambda)+s_0\xi^2u_\lambda-(g_0+\lambda)e^{\frac{2}{n}u_\lambda}\xi^2u_\lambda\right)d\mu_\eta=0,
\end{align}
here $\xi\in C^\infty_0(M)$ satisfying 
\begin{align}
\xi(x)=
\begin{cases}
1& x\in D\\
0& x\in M\setminus D'
\end{cases}
,\ 0\leq\xi\leq1,\ |\nabla \xi|\leq A,
\end{align}
where $D\subset\subset D'\subset\subset M_-$, $A$ is a constant related to $D, D'$. Note that since $xe^x\geq-1$ for any $x\in \mathbb{R}$ and $D'\subset\subset M_-$, it follows that
\begin{align}
\int_{D'} (g_0+\lambda)e^{\frac{2}{n}u_\lambda}\xi^2u_\lambda d\mu_\eta\leq -\frac{n}{2} \inf_{D'}g_0\text{Vol}(D').
\end{align}
Then, by a similar argument as in Lemma \ref{lemma4.6}, it is easy to get
\begin{align}\label{3.71..}
\int_D |\nabla u_\lambda|^2d\mu_\eta&\leq \int_{D'} |\nabla (\xi u_\lambda)|^2d\mu_\eta\\
&=\int_{D'}\left(|\nabla \xi|^2|u_\lambda|^2-s_0\xi^2u_\lambda+(g_0+\lambda)e^{\frac{2}{n}u_\lambda}\xi^2u_\lambda\right)d\mu_\eta\notag\\
&\leq A^2 \|u_\lambda\|_{L^2(D')}^2+ \|s_0\|_{L^2(D')}\|u_\lambda\|_{L^2(D')}-\frac{n}{2} \inf_{D'}g_0\text{Vol}(D').\notag
\end{align}
Combining \eqref{3.67..} and \eqref{3.71..} yields that $\|u_\lambda\|_{W^{1,2}(D)}$ is uniformly bounded as $\lambda\rightarrow \lambda^*$.
\qed

\begin{lemma}
$u_\lambda$ is locally uniformly $C^0$-bounded in the open set $M_-$ as $\lambda\rightarrow \lambda^*$.
\end{lemma}
\proof We take a decreasing sequence $\{r_i\}_{i=0}^\infty$ in $\mathbb{R}^+$, and choose the following cut-off functions $\xi\in C^\infty_0(M)$ at each step:
\begin{align}\label{4.38..}
\xi(x)=
\begin{cases}
1& x\in B_{r_{i+1}}(x_0)\\
0& x\in M\setminus B_{r_ i}(x_0)
\end{cases}
,\ 0\leq\xi\leq1,\ |\nabla \xi|\leq \frac{C}{r_i-r_{i+1}},
\end{align}
where $B_{r}(x_0)=\{x\in M\ :\ dist(x_0, x)<r\}$, $r_i=\theta+\frac{\tau-\theta}{2^i}$, $\theta<r_i\leq\tau$ for $i=0,1,2, \cdots$, and $B_\tau(x_0)\subset M_-$. Pick the test function $\varphi=\xi^2u_\lambda|u_\lambda|^{a}\in C^1_0(M)\ (a\geq0)$ in 
\begin{align}
\int_M \left(\nabla u_\lambda\cdot \nabla \varphi+ s_0\varphi- (g_0+\lambda)e^{\frac{2}{n}u_\lambda}\varphi\right)d\mu_\eta=0.\label{3.68}
\end{align}
By a simple computation, we obtain
\begin{align}\label{3.69}
|\nabla(\xi u_\lambda|u_\lambda|^{\frac{a}{2}})|^2=&\left(\frac{a}{2}+1\right)\nabla u_\lambda\cdot \nabla (\xi^2u_\lambda|u_\lambda|^{a})+|\nabla\xi|^2|u_\lambda|^{a+2}\\
&\ \ \ \ \ \ -\frac{a}{2}\left(\frac{a}{2}+1\right)\xi^2| u_\lambda|^a|\nabla u_\lambda|^2.\notag
\end{align}
In terms of \eqref{3.68} and \eqref{3.69}, we have
\begin{align}
\int_M& |\nabla(\xi u_\lambda|u_\lambda|^{\frac{a}{2}})|^2d\mu_\eta\\
&\leq \left(\frac{a}{2}+1\right)\int_M \nabla u_\lambda\cdot \nabla (\xi^2u_\lambda|u_\lambda|^{a})d\mu_\eta+\int_M |\nabla\xi|^2|u_\lambda|^{a+2}d\mu_\eta \notag\\
&\leq \left(\frac{a}{2}+1\right)\left\{\int_{B_{r_i}}(g_0+\lambda)e^{\frac{2}{n}u_\lambda}\xi^2u_\lambda|u_\lambda|^ad\mu_\eta-\int_{B_{r_i}}s_0\xi^2u_\lambda|u_\lambda|^ad\mu_\eta\right\}\notag\\
&\ \ \ \ \ \ \ \ \ +\frac{C^2}{(r_i-r_{i+1})^2}\int_{B_{r_i}}|u_\lambda|^{a+2}d\mu_\eta.\notag
\end{align}
Since $B_{r_i}\subset M_-$ and the fact that $xe^x\geq-1$ for any $x\in \mathbb{R}$, we deduce that
\begin{align}
\int_{B_{r_i}}(g_0+\lambda)e^{\frac{2}{n}u_\lambda}\xi^2u_\lambda|u_\lambda|^ad\mu_\eta\leq C_1\int_{B_{r_i}}|u_\lambda|^ad\mu_\eta,
\end{align}
where $C_1= \frac{n}{2}\max_M{(-g_0)}$. Besides, we have the following
\begin{align}
-\int_{B_{r_i}}s_0\xi^2u_\lambda|u_\lambda|^ad\mu_\eta\leq C_2\int_{B_{r_i}}|u_\lambda|^{a+1}d\mu_\eta,
\end{align}
where $C_2=\max_M |s_0|$. Hence, 
\begin{align}\label{3.73}
\int_M& |\nabla(\xi u_\lambda|u_\lambda|^{\frac{a}{2}})|^2d\mu_\eta\\
&\leq \left(\frac{a}{2}+1\right)\left\{C_1\int_{B_{r_i}}|u_\lambda|^ad\mu_\eta+C_2\int_{B_{r_i}}|u_\lambda|^{a+1}d\mu_\eta \right\}\notag\\
&\ \ \ \ \ \ \ \ \ +\frac{C^2}{(r_i-r_{i+1})^2}\int_{B_{r_i}}|u_\lambda|^{a+2}d\mu_\eta\notag\\
&\leq C_3\left(a+\frac{1}{2}+\frac{1}{(r_i-r_{i+1})^2} \right) \int_{B_{r_i}}|u_\lambda|^{a+2}d\mu_\eta+\frac{3}{2}C_3,\notag
\end{align}
where we have used the inequalities \eqref{3.16...} and \eqref{3.17...}, and $C_3=\max\{C_1, C_2, C^2\}$. Applying the Sobolev’s inequality (cf. Lemma \ref{theorem2.3}) to $f=\xi u_\lambda|u_\lambda|^{\frac{a}{2}}$ yields
\begin{align}\label{3.74}
\left(\int_M(\xi u_\lambda|u_\lambda|^{\frac{a}{2}})^{2\beta}d\mu_\eta\right)^{\frac{1}{\beta}}\leq C_4\int_M \left(|\nabla(\xi u_\lambda|u_\lambda|^{\frac{a}{2}})|^2+\xi^2|u_\lambda|^{a+2}\right)d\mu_\eta,
\end{align}
where $\beta=\frac{n}{n-1}$. Combining \eqref{3.73} and \eqref{3.74}, we have
\begin{align}\label{3.75.}
&\left(\int_{B_{r_{i+1}}}|u_\lambda|^{(a+2)\beta}d\mu_\eta\right)^{\frac{1}{\beta}}\\
&\ \ \ \ \ \ \ \leq C_5\left(a+2+\frac{1}{(r_i-r_{i+1})^2} \right)\int_{B_{r_i}}|u_\lambda|^{a+2}d\mu_\eta+\frac{3}{2}C_5\notag\\
&\ \ \ \ \ \ \ \leq C_6\left(a+2+\frac{1}{(r_i-r_{i+1})^2} \right)\max\left\{\int_{B_{r_i}}|u_\lambda|^{a+2}d\mu_\eta,1\right\},\notag
\end{align}
where $C_5=\max\{C_4C_3, C_4\}$, $C_6=2C_5$ and $a\geq 0$. Let us choose $a_i$ in \eqref{3.75.} such that $\beta^i=\frac{a_i}{2}+1$, then
\begin{align}
\max&\{\|u_\lambda\|_{L^{2\beta^{i+1}}(B_{r_{i+1}})},1\}\\
&\leq C_6^{\frac{1}{2\beta^i}}\left( 2\beta^i+\frac{4^{i+1}}{(\tau-\theta)^2}\right)^{\frac{1}{2\beta^i}}\max\{\|u_\lambda\|_{L^{2\beta^i}(B_{r_i})},1\},\notag\\
&\leq C_6^{\frac{1}{2\beta^i}}4^{\frac{i+1}{2\beta^i}}\left(\frac{1}{2}+\frac{1}{(\tau-\theta)^2}\right)^{\frac{1}{2\beta^i}}\max\{\|u_\lambda\|_{L^{2\beta^i}(B_{r_i})},1\}, \notag
\end{align}
for any $i=0,1,2,\cdots$, since $\frac{\beta}{4}=\frac{n}{4(n-1)}<1$. By iteration, we conclude that
\begin{align}
\max&\{\|u_\lambda\|_{L^{2\beta^{i+1}}(B_{r_{i+1}})},1\}\\
&\leq C_6^{\frac{1}{2\beta^i}}4^{\frac{i+1}{2\beta^i}}C_7^{\frac{1}{2\beta^i}}\max\{\|u_\lambda\|_{L^{2\beta^i}(B_{r_i})},1\}\notag\\
&\leq\dots\notag\\
&\leq C_6^{\frac{1}{2}\sum_{k=0}^i\frac{1}{\beta^k}}4^{\frac{1}{2}\sum_{k=0}^i\frac{k+1}{\beta^k}}C_7^{\frac{1}{2}\sum_{k=0}^i\frac{1}{\beta^k}}\max\{\|u_\lambda\|_{L^{2}(B_{\tau})},1\}  \notag
\end{align}
where $C_7=\frac{1}{2}+\frac{1}{(\tau-\theta)^2}$. Now let $i\rightarrow +\infty$, we obtain
\begin{align}\label{3.78}
\max\{\|u_\lambda\|_{C^{0}(B_{\theta})},1\}\leq C_8\max\{\|u_\lambda\|_{L^{2}(B_{\tau})},1\}\leq C_9, 
\end{align}
where $C_8, C_9$ are two constants independent of $\lambda$, and the second inequality is due to Lemma \ref{lemma4.7}.
\qed

\begin{lemma}\label{lemma4.9}
$e^{\frac{u_\lambda}{n}}$ is uniformly $W^{1,2}(M)$-bounded as $\lambda\rightarrow \lambda^*$.
\end{lemma}
\proof Let $h$ be a nonpositive smooth function on $M$ satisfying $h<0$ in some open subset $D$ of $M_-$, $h\equiv 0$ in $M\setminus D$. Since $h\leq (\not\equiv) 0$, then there is a unique solution $v\in C^\infty(M)$ of
\begin{align}\label{4.49}
-\Delta^{Ch}_\eta v+s_0=he^{\frac{2}{n}v}.
\end{align}
Set $\omega_\lambda=u_\lambda-v$, then from \eqref{3.49.} and \eqref{4.49} it satisfies
\begin{align}\label{3.85}
-\Delta^{Ch}_\eta \omega_\lambda=(g_0+\lambda)e^{\frac{2}{n}u_\lambda}-he^{\frac{2}{n}v}.
\end{align}
Multiplying \eqref{3.85} by $e^{\frac{2}{n}\omega_\lambda}$ and integrating on $M$ yield
\begin{align}
2n\int_M |\nabla e^{\frac{\omega_\lambda}{n}}|^2d\mu_\eta-\int_M (g_0+\lambda)e^{\frac{2}{n}(u_\lambda+\omega_\lambda)}d\mu_\eta=-\int_M h e^{\frac{2}{n}u_\lambda}d\mu_\eta,
\end{align}
and then using \eqref{3.52.} with $\varphi=e^{\frac{1}{n}\omega_\lambda}$, we conclude 
\begin{align}\label{4.53..}
\frac{3n}{2}\int_M  |\nabla e^{\frac{\omega_\lambda}{n}}|^2d\mu_\eta\leq -\int_Dh e^{\frac{2}{n}u_\lambda}d\mu_\eta\leq -(\inf_D h) e^{\frac{2}{n}\|u_\lambda\|_{L^\infty(D)}}d\mu_\eta\leq C_{10},
\end{align}
where $C_{10}$ is a constant independent of $\lambda$, and last inequality holds because $u_\lambda$ is locally uniformly $L^\infty$-bounded in the open set $M_-$. Now we claim that $e^{\frac{\omega_\lambda}{n}}$ is uniformly bounded in $L^2(M)$. Indeed, if it is not true, then there is a subsequence $\lambda_i\rightarrow \lambda^*$ such that $\lim_{i\rightarrow +\infty}\|e^{\frac{\omega_{\lambda_i}}{n}}\|_{L^2(M)}=+\infty$. Set 
\begin{align}
v_{\lambda_i}=\frac{e^{\frac{\omega_{\lambda_i}}{n}}}{\|e^{\frac{\omega_{\lambda_i}}{n}}\|_{L^2(M)}},
\end{align}
then $\|v_{\lambda_i}\|_{L^2(M)}=1$. By \eqref{4.53..}, we obtain
\begin{align}
\lim_{i\rightarrow +\infty}\|\nabla v_{\lambda_i}\|_{L^2(M)}=0.
\end{align}
Hence, passing to a subsequence, $v_{\lambda_i}\rightharpoonup v_{\lambda^*}$ weakly in $W^{1,2}(M)$, where $v_{\lambda^*}\equiv C^*$ is a constant with $\|v_{\lambda^*}\|_{L^2(M)}=1$. On the other hand, since $\|u_\lambda\|_{L^\infty(D)}$ is uniformly bounded as $\lambda\rightarrow \lambda^*$, it follows that $\lim_{i\rightarrow +\infty}\|v_{\lambda_i}\|_{L^2(D)}= 0$, thus $C^*\equiv 0$, which leads to a contradiction with $\|v_{\lambda^*}\|_{L^2(M)}=1$. According to the relationship $e^{\frac{u_\lambda}{n}}=e^{\frac{v}{n}}e^{\frac{\omega_\lambda}{n}}$, we get $\|e^{\frac{u_\lambda}{n}}\|_{W^{1,2}(M)}$ is also uniformly bounded as $\lambda\rightarrow \lambda^*$.
\qed

\begin{lemma}\label{lemma4.10}
If $\dim_{\mathbb{C}} M=n=2$, then $u_\lambda$ is uniformly $C^{2,\alpha}(M)$-bounded as $\lambda\rightarrow \lambda^*$.
\end{lemma}
\proof By the Sobolev embedding theorem $W^{1,2}(M)\subset L^{4}(M)$ ($\dim_{\mathbb{C}}M=2$) and Lemma \ref{lemma4.9}, we have
\begin{align}
\int_M e^{2u_\lambda}d\mu_\eta\leq C_{11},
\end{align}
where $C_{11}$ is a constant independent of $\lambda$. Hence, 
\begin{align}\label{4.56}
\|(g_0+\lambda)e^{u_\lambda}\|_{L^2(M)}\leq (\max_M|g_0|+\lambda^*)\|e^{u_\lambda}\|_{L^2(M)}\leq (\max_M|g_0|+\lambda^*)C_{11}.
\end{align}
On the other hand,
\begin{align}\label{4.57}
\int_M |u_\lambda|^2d\mu_\eta&\leq \int_{M\cap \{-C\leq u_\lambda\leq0\}}|u_\lambda|^2d\mu_\eta+\int_{M\cap \{u_\lambda>0\}}|u_\lambda|^2d\mu_\eta\\
&\leq C^2\text{Vol}(M)+\int_M e^{2u_\lambda}d\mu_\eta\leq C^2+C_{11},\notag
\end{align}
where $C$ is a constant defined as in Lemma \ref{lemma4.4}. From \eqref{4.56}, \eqref{4.57} and the elliptic $L^2$-estimate for \eqref{3.49.}, it follows that
\begin{align}\label{4.58}
\|u_\lambda\|_{W^{2,2}(M)}\leq C\left(\|u_\lambda\|_{L^2(M)} +\|(g_0+\lambda)e^{u_\lambda}\|_{L^2(M)}\right)\leq C_{12},
\end{align}
where $C_{12}$ is a constant independent of $\lambda$. Now we recall the Trudinger-Moser inequality (cf. Lemma \ref{lemma2.4}): there exist two positive constants $\delta_1$ and $\delta_2$ such that
\begin{align}\label{4.59}
\int_Me^{\delta_1 \left(\frac{u-\bar{u}}{\|\nabla^2u\|_{L^2(M)}}\right)^2}d\mu_\eta\leq \delta_2
\end{align}
for any $u\in W^{2,2}(M)$, where $\bar{u}=\int_M ud\mu_\eta$. Then, using\eqref{4.57}, \eqref{4.58} and \eqref{4.59} yields
\begin{align}\label{3.95}
\int_Me^{pu_\lambda}d\mu_\eta&\leq e^{p\bar{u}_\lambda}e^{\frac{p^2}{4\delta_1}\|\nabla^2 u_\lambda\|_{L^2(M)}^2}\int_Me^{\delta_1 \left(\frac{u_\lambda-\bar{u}_{\lambda}}{\|\nabla^2u_{\lambda}\|_{L^2(M)}}\right)^2}d\mu_\eta\\
&\leq e^{p\sqrt{C^2+C_{11}}}e^{\frac{p^2}{4\delta_1}C_{12}^2}\delta_2\notag
\end{align}
 for any $p>0$. According to \eqref{3.95} and the elliptic $L^p$-estimate, we see that $u_\lambda$ is uniformly $W^{2,p}(M)$-bounded as $\lambda\rightarrow \lambda^*$. If we pick $p>4=\dim_{\mathbb{R}} M$, then $u_\lambda$ is uniformly $C^{1,\alpha}(M)$-bounded as $\lambda\rightarrow \lambda^*$. Using the Schauder estimate for \eqref{3.49.}, we have $u_\lambda$ is uniformly $C^{2,\alpha}(M)$-bounded as $\lambda\rightarrow \lambda^*$.
 \qed
 
\proof[\textbf{Proof of Theorem \ref{theorem4.2}}] According to Lemma \ref{lemma4.10} and the compactly embedding theorem $C^{2,\alpha}(M)\subset\subset C^{2,\beta}(M)\  (0<\beta<\alpha<1)$, we obtain that there is a subsequence $\lambda_i\rightarrow \lambda^*$ such that $u_{\lambda_i}$ converges to $u_{\lambda^*}$ in $C^{2,\beta}(M)$ for any $\beta\in (0,1)$. Hence, we find a solution $u_{\lambda^*}$ for \eqref{3.49.} with $\lambda=\lambda^*$, and by the regularity results for elliptic equations, we know that $u_{\lambda^*}$ is smooth on $M$.
\qed

\bigskip
\bigskip

Weike Yu

School of Mathematics and Statistics

Nanjing University of Science and Technology

Nanjing, 210094, Jiangsu, P. R. China

wkyu2018@outlook.com

\bigskip

\end{document}